\documentclass[12pt,a4paper]{article}
\usepackage{amssymb}
\usepackage{amsfonts}
\usepackage{srcltx}
\usepackage{amsmath}
\usepackage{amsthm}
\usepackage{enumerate}
\usepackage{indentfirst}

\newtheorem{theorem}{Theorem}[section]

\newtheorem{lemma}[theorem]{Lemma}

\numberwithin{equation}{section}

\begin{document}
\title{
On an extremal problem for nonoverlapping domains
\thanks{The authors is supported by
the Mendeleev Grant of National Research Tomsk State University
and RFBR Grant No 16-01-00121 À
}}
\author{Pchelintsev E.A.
\thanks{
Department of Mathematics and Mechanics,
National Research Tomsk State University,
Lenin str. 36,
 634050 Tomsk, Russia,
 e-mail: evgen-pch@yandex.ru
}
 \and
Pchelintsev V.A.\thanks{
 Department of Higher Mathematics and Mathematical Physics,
 Tomsk Polytechnic University,
Lenin str. 30,
 634050 Tomsk, Russia,
 e-mail: vpchelintsev@vtomske.ru}
}
 \date{}
\maketitle

\begin{abstract}
The paper considers the problem of finding the range of functional
$I=J\left(f(z_0),\overline{f(z_0)},
F(\zeta_0),\overline{F(\zeta_0)}\right)$,
defined on the class $\mathfrak{M}$ of pairs functions
 $(f(z),F(\zeta))$ that are univalent in the system of the disk and the interior of the disk, using the method
 of internal variations.
We establish that the range of this functional is bounded by
the curve whose equation is written in terms of elliptic integrals,
depending on the parameters of the functional $I$.
\end{abstract}

{\bf Key words:}  Method of internal variations, Univalent function, Nonoverlapping domains, Functional range, Elliptic integrals
\\
\par
{\bf AMS (1991) Subject Classification : 30C70}

\bibliographystyle{plain}

%\newpage

\section{Introduction}\label{sec:In}

In the geometric theory of univalent functions there are a number of papers
and books devoted to the the problem of nonoverlapping domains.
These problems were developed by M.A. Lavrentiev~\cite{Lavrentiev}
G.M. Goluzin~\cite{Goluzin}, James A. Jenkins~\cite{Jenkins}, Z. Nehari~\cite{Nehari},
N.A. Lebedev~\cite{Lebedev1957}, M. Shiffer~\cite{Duren},
R. K\"{u}hnau~\cite{Kuhnau} and others.
Summary of results in this area contained in~\cite{Bakhtin}.

Let $D$ and $D^{*}$ be nonoverlapping simply connected domains in the $w$-plane such that
 $0\in D$ and $\infty\in D^{*}$.
Assume that $f:E\to D$ and $F:E^*\to D^{*}$ are
holomorphic and meromorphic respectively univalent functions
normalized by the conditions
 $f(0)=0$ and
$F(\infty)=\infty$. Here
$E=\{z \in \mathbb{C}:|z|<1\}$ and
$E^{*}=\{\zeta \in \mathbb{C}:|\zeta|>1\}$.
The family of all such pairs
$(f(z),F(\zeta))$ is called the class
$\mathfrak{M}$. Some extremal problems on this class were studied in
\cite{Lebedev1957, Andreev, Pchelintsev}.

Let $J:G\rightarrow \mathbb{C}$,
$J=J(\omega_1,\omega_2,\omega_3,\omega_4)$ is an
analytic in some domain $G \subset \mathbb{C}^{4}$ nonhomogeneous and nonconstant function.

Now we fix an arbitrary point
$z_{0}\in E$, $\zeta_{0}\in E^{*}$ and
define on the class $\mathfrak{M}$
a functional
\begin{equation}\label{eq2.1}
I:\mathfrak{M} \rightarrow \mathbb{C},\,
\quad
I(f,F)=J\left(f(z_0),\overline{f(z_0)},
F(\zeta_0),\overline{F(\zeta_0)}\right).
\end{equation}

This paper considers the problem of finding the range
$\Delta$ of the functional $I$ on the class
$\mathfrak{M}$. Since together with pair $(f(z),F(\zeta))$ the class
$\mathfrak{M}$ contains the following pairs of functions
$(f(z e^{i\varphi}),F(\zeta e^{i\psi}))$
for any parameters
$\varphi,\,\psi\in\mathbb{R}$, then we reduce the initial problem
to equivalent one for the functional
$$
I=J\left(f(r),\overline{f(r)},
F(\rho),\overline{F(\rho)}\right),
$$
where $r=|z_{0}|\in(0,1)$,
$\rho=|\zeta_{0}|\in(1,+\infty)$.
Further we will solve this problem.

Taking into account that the class $\mathfrak{M}$ contains the pair of function
$(f(zt)$, $F(\zeta t^{-1}))$ for any $t \in [0, 1]$ (see, \cite{Andreev}),
we have that the range $\Delta$
of the functional \eqref{eq2.1} is connected set.
Note that if the class $\mathfrak{M}$ complement the pair of functions
$(f(z), \infty)$, where $f(z)$, $f(0)=0$, is holomorphic univalent
function in $E$, then $\Delta$ will be
closed set \cite{Ulina}. Hence,
it suffices to find the boundary $\Gamma$
of the set $\Delta$.

A point $I_{0}\in \Gamma$ is called a \textit{nonsingular boundary point} if there exists an exterior point $I_{e}$ for $\Delta$
such that the distance between $I_{e}$ and $I_{0}$  is equal to the distance between $I_{e}$ and the set
$\Delta$, i.e.
\begin{equation}\label{eq2.2}
|I_{0}-I_{e}|=
\inf \limits_{I \in \Delta}|I-I_{e}|.
\end{equation}

The set $\Gamma_{0}$ of nonsingular boundary points
is dense in
$\Gamma$~\cite{Lebedev1955}. Thus,
the initial problem of finding the range $\Delta$ of \eqref{eq2.1}
 is replaced by an equivalent extremal problem:
find the minimum of the real-valued functional
$|I-I_{e}|$ on $\mathfrak{M}$  for all possible $I_{e} \notin\Delta$.

The functions giving nonsingular boundary points of a functional are called the \textit{boundary functions} of
this functional; i.e., these are the functions at which the values of the functional are nonsingular boundary
points.

In this paper to solve the problem we apply the Schiffer's method of internal variations \cite{Schiffer} using pairs of varied functions from \cite{Lebedev1957}.

Main results are given in the following Section 2.

\section{Differential equations for the boundary functions.
The equation of the boundary of $\Delta$}

Write down the variational formulas for boundary functions $f(z)$ and $F(\zeta)$
on the class $\mathfrak{M}$ in the form
$$
f_{\varepsilon}(z)=f(z)+\varepsilon P(z)+
o(z,\varepsilon),
$$
$$
F_{\varepsilon}(\zeta)=F(\zeta)+
\varepsilon Q(\zeta)+
o(\zeta,\varepsilon)
$$
for $\varepsilon$ positive and sufficiently small.
In view of the functional \eqref{eq2.1} is
G\^{a}teaux differentiable, then we can it rewrite in the form
\begin{equation*}
I_{*}=I+
\varepsilon
\left\{\frac{\partial J(\omega_0)}
{\partial \omega_1}P(r)+
\frac{\partial J(\omega_0)}
{\partial \omega_2}\overline{P(r)}+
\frac{\partial J(\omega_0)}
{\partial \omega_3}Q(\rho)+
\frac{\partial J(\omega_0)}
{\partial \omega_4}\overline{Q(\rho)}\right\}+
o(\varepsilon),
\end{equation*}
where $I_{*}=
I(f_{\varepsilon}(z),F_{\varepsilon}(\zeta))$,
$\omega_0=\left(f(r),\overline{f(r)},
F(\rho),\overline{F(\rho)}\right) \in G$.

Let $(f(z),F(\zeta))$ be a boundary pair of functions of this functional giving the point $I_{0}$.
The equality \eqref{eq2.2} implies that
$$
|I_{*}-I_{e}|\geq |I_{0}-I_{e}|.
$$
From here, hence that each boundary pair of functions satisfies the necessary condition
\begin{equation}\label{eq3.1}
\mathrm{Re}
[pP(r)+qQ(\rho)]\geq0,
\end{equation}
where
$$
p=e^{-i\alpha}\frac{\partial J(\omega_0)}
{\partial \omega_1}+e^{i\alpha}
\overline{\left(\frac{\partial J(\omega_0)}
{\partial \omega_2}\right)}, \quad
q=e^{-i\alpha}\frac{\partial J(\omega_0)}
{\partial \omega_3}+e^{i\alpha}
\overline{\left(\frac{\partial J(\omega_0)}
{\partial \omega_4}\right)},
$$
and $\alpha=\arg(I-I_{e}).$

\begin{lemma}
Let $f(z)$, $F(\zeta)$ be boundary functions of functional \eqref{eq2.1}.
Then the union of domains $D$
and $D^{*}$ has no exterior points in $\overline{\mathbb{C}}_w$.
\end{lemma}

{\bf Proof.}
Suppose that $D\cup D^{*}$ has at least one exterior point $w_0$ in $w$-plane.
Then by definition of the exterior point there exists
neighborhood of the point $w_0$
consisting of exterior points. Using the varied function
\begin{equation*}
f(z,\varepsilon)=f(z)+\varepsilon
A_0\frac{f(z)}{f(z)-w_{0}}\,,\;
F(\zeta,\varepsilon)=F(\zeta)+\varepsilon
A_0\frac{F(\zeta)}{F(\zeta)-w_{0}}\,;
\end{equation*}
where $w_0$ is an exterior point for $D$ and $D^{*}$ simultaneously and
$A_0$ is an arbitrary complex constant, as the comparison function
in \eqref{eq3.1}, rewrite it as
\begin{equation}\label{eq3.2}
\mathrm{Re}\left(A_0
\frac{R(w_0)}
{(f(r)-w_0)(F(\rho)-w_0)}\right)\geq0,
\end{equation}
where $R(w_0)$ is a  linear polynomial.
The fraction in the condition \eqref{eq3.2} is equal to zero,
otherwise in view of the arbitrariness of $\arg A_0$
we can choose it so that the left side
will be negative. Therefore $R(w_0)=0$.
It is possible only for a single point, while inequality
\eqref{eq3.2} should be performed for any point from neighborhood of the $w_0$.
This contradiction proves the lemma.

Further, to obtain differential equations for the boundary functions of the functional \eqref{eq2.1}
we consider the following pairs of variational formulas from \cite{Lebedev1957}:

\noindent 1)
\begin{equation*}
f(z,\varepsilon)=f(z)+\varepsilon
A_0\left(\frac{f(z)}{f(z)-f(z_{0})}-
\frac{f(z_0)}{z_0f'^{2}(z_0)}\frac{z
f'(z)}{z-z_0}\right)+
\end{equation*}

\begin{equation}\label{3.2.1}
+\varepsilon\overline{A_0}
\overline{\frac{f(z_{0})}{z_0f'^{2}(z_0)}}
\frac{z^2f'(z)}{1-\overline{z_0}z}+
o(z,\varepsilon),
\end{equation}

\begin{equation*}
F(\zeta,\varepsilon)=F(\zeta)+\varepsilon
A_0\frac{F(\zeta)}{F(\zeta)-f(z_0)}\,,
\end{equation*}
where $z_0\in E$,\ $A_0$\ is an arbitrary complex
constant;

\noindent 2)
\begin{equation*}
f(z,\varepsilon)=f(z)+\varepsilon
A_0\frac{f(z)}{f(z)-F(\zeta_0)},
\end{equation*}

\begin{equation}\label{3.2.2}
F(\zeta,\varepsilon)=F(\zeta)+
\varepsilon A_0\left(
\frac{F(\zeta)}{F(\zeta)-F(\zeta_0)}-
\frac{F(\zeta_0)}{\zeta_0^2F'^2(\zeta_0)}
\frac{\zeta^2F'(\zeta)}
{\zeta-\zeta_0}\right)+
\end{equation}

\begin{equation*}
+\varepsilon\overline{A_0}
\overline{\frac{F(\zeta_0)}
{\zeta_0^2F'^2(\zeta_0)}}
\frac{\zeta
F'(\zeta)}{1-\overline{\zeta_0}\zeta}+
o(\zeta,\varepsilon),
\end{equation*}
where $\zeta_0\in E^*$, $A_0$ is an arbitrary complex
constant.

\begin{theorem}\label{Th.1}
Every boundary pair of functions
$(f(z),F(\zeta))$ of the functional
\eqref{eq2.1} satisfies in
$E$ and $E^{*}$ the system
of functional-differential equations
\begin{equation}\label{eq3.3}
\frac{(C_1f(z)-C_2)f'(z)}
{f(z)(f(z)-f(r))(f(z)-F(\rho))}=
\frac{A}{z(r-z)(1-rz)}\,,
\end{equation}
\begin{equation}\label{eq3.4}
\frac{(C_1F(\zeta)-C_2)F'(\zeta)}
{F(\zeta)(F(\zeta)-f(r))(F(\zeta)-F(\rho))}=
\frac{B}{\zeta(\rho-\rho)(1-\rho \zeta)}\,,
\end{equation}
where
$$
C_1=pf(r)+qF(\rho), \quad
C_2=(p+q)f(r)F(\rho),
$$
$$
A=(1-r^2)rpf'(r)>0, \quad
B=(\rho^{2}-1) \rho F'(\rho)>0.
$$

\end{theorem}

{\bf Proof.}
If we choose the variational formula \eqref{3.2.1}, then
\eqref{eq3.1} takes the form
\begin{multline*}
\mathrm{Re}\biggl[p A_0\frac{
 f(r)}{f(r)-f(z_{0})}- p A_0
 \frac{r f'(r)}{r-z_0}
\frac{f(z_0)}{z_0f'^{2}(z_0)}+
\\[3mm]
+p \overline{A_0}
\frac{r^{2}f'(r)}{1-r
\overline{z_0}}\overline{\frac{f(z_{0})}
{z_0f'^{2}(z_0)}}+q A_0\frac{
F(\rho)}{F(\rho)-f(z_0)}\biggr]\geq0.
\end{multline*}
Replacing the last summand under the real part by its conjugate, we have
\begin{multline*}
\mathrm{Re}A_0\biggl[\frac{p
 f(r)}{f(r)-f(z_{0})}-
 p\frac{r f'(r)}{r-z_0}
\frac{f(z_0)}{z_0f'^{2}(z_0)}+
\\[3mm]
+ \overline{p}\overline{\frac{r^{2}f'(r)}{1-r
\overline{z_0}}}\frac{f(z_{0})}
{z_0f'^{2}(z_0)}+\frac{q
F(\rho)}{F(\rho)-f(z_0)}\biggr]\geq0.
\end{multline*}
In this condition, the expression in parentheses is equal to zero; otherwise, under an appropriate choice
of $\arg A_0$ we would get that the left-hand side of the last inequality is negative. This leads to the equality
$$
\frac{pf(r)}{f(r)-f(z_{0})}+
\frac{qF(\rho)}{F(\rho)-f(z_0)}=
\frac{f(z_0)}{z_0f'^{2}(z_0)}
\left(\frac{rp f'(r)}{r-z_0}-
\overline{p f'(r)}
\frac{r^{2}}{1-rz_0} \right).
$$
Since, in this relation,
$z_0$  is an arbitrary point of
$E$, replacing $z_0$ by $z$, and  in view of $pf'(r)<0$, we obtain
a differential equation
for the boundary function $f(z)$.
The calculations show that it has the form \eqref{eq3.3}.

Ideologically, the deduction of \eqref{eq3.4}, repeats that \eqref{eq3.3};
 for this we must apply \eqref{eq3.1}  together with the
variational formulas \eqref{3.2.2} and use inequality $qF'(\rho)>0$.
The theorem is proved.

From the analytic theory of differential equations \cite{Golubev}, we conclude that the boundary functions
$f(z)$ and $F(\zeta)$ satisfying their equations are holomorphic not only in
 $E$ and $E^{*}$,  but also on the unit circle $|z|=|\zeta|=1$. From here and in view that
the union
$D\cup D^{*}$ do not
contain exterior points, one has that the domains
$D$ and $D^{*}$
are bounded by some closed analytic Jordan curve.

Further, to find the equation of the boundary of the range
$\Delta$ of the functional \eqref{eq2.1} we integrate
\eqref{eq3.3} and \eqref{eq3.4}.

Extract the square root from both sides of \eqref{eq3.3} and integrate the result by
$z$ from $0$ to $r$.
Consider the left-hand side :
$$
J=\int\limits_0^r
\sqrt{\frac{C_1f(z)-C_2}{f(z)(f(z)-f(r))(f(z)-F(\rho))}}
\,f'(z)dz.
$$
Changing the integration variable $t=f(z)/f(r)$, we have
$$
J=a\int\limits_0^1
\frac{t-b}{\sqrt{t(1-t)(1-ct)(t-b)}}\,dt,
$$
where
$$
a=\sqrt{C_1 \frac{f(r)}{F(\rho)}}, \quad
b=\frac{C_2}{C_1f(r)}, \quad
c=\frac{f(r)}{F(\rho)}.
$$
Putting $t=1/x$ in $J$, we infer
$$
J=a\int\limits_1^\infty
\frac{dx}{x\sqrt{(x-1)(x-c)(1-bx)}}-ab
\int\limits_1^\infty
\frac{dx}{\sqrt{(x-1)(x-c)(1-bx)}}\,.
$$
Performing the change of variables $y=b(x-1)/(bx-1)$ in the integrals,
after transformations one obtains
$$
J=-2ph \mathbf{\Pi}(n,k),
$$
where
$$
\mathbf{\Pi}(n,k)=\int\limits_0^
{\pi/2}\frac{dt}{(1+n\sin^2t)\sqrt{1-k^2\sin^2t}}
$$
is the complete elliptic integral of the third kind,
$$
h=\sqrt{\frac{(f(r)-
F(\rho))f^2(r)}{C_2}}\,,\quad
n=-\frac{1}{b}\,,\quad
k=\sqrt{\frac{q}{p+q}}\,.
$$
Here $\sqrt{1-k^2\sin^2t}$ stands for the branch of the function assuming $1$ at $t \to 0$.

Now integrate the right-hand side:
$$
J=\sqrt{A}\int\limits_0^r
\frac{dz}{\sqrt{z(r-z)(1-rz)}}.
$$
Changing the integration variable $x=z/r$, one has
$$
J=2\sqrt{A}\,\mathbf{K}(r),
$$
where
$$
\mathbf{K}(r)=\int\limits_0^{\pi/2}
\frac{dt}{\sqrt{1-r^2\sin^2t}}
$$
is the complete elliptic integral of the first kind.

Thus, upon integration,
\eqref{eq3.3}
looks as
\begin{equation}\label{eq3.5}
-ph \mathbf{\Pi}(n,k)=\sqrt{A}\,\mathbf{K}(r).
\end{equation}

Integrate \eqref{eq3.4}
after extracting the square root with respect to
$\zeta$ from $\rho$ to $\infty$.
First consider the left-hand side:
$$
L=\int\limits_{\rho}^\infty
\sqrt{\frac{C_1F(\zeta)-C_2}
{F(\zeta)(F(\zeta)-F(\rho))(F(\zeta)-f(r))}}
\,F'(\zeta)\,d\zeta.
$$
Changing the integration variable
$t=F(\zeta)/F(\rho)$, we have
$$
L=a^{*}\int\limits_1^\infty
\frac{t-b^{*}}{\sqrt{t(1-t)(1-c^{*}t)(t-b^{*})}}\,dt,
$$
where
$$
a^{*}=\sqrt{C_1 \frac{F(\rho)}{f(r)}}, \quad
b^{*}=\frac{C_2}{C_1 F(\rho)},
\quad c^{*}=\frac{F(\rho)}{f(r)}.
$$
Putting $t=1/x$ in $L$, we infer
$$
L=a^*\int\limits_0^1
\frac{dx}{x\sqrt{(x-1)(x-c^{*})(1-b^{*}x)}}-
a^*b^*\int\limits_0^1
\frac{dx}{\sqrt{(x-1)(x-c^{*})(1-b^{*}x)}}\,.
$$
Performing the change of variables $u=b^{*}(1-x)/(b^{*}-1)$ in the integrals,
after calculation one comes to equality
$$
L=2\left(l\mathbf{\Pi}(\phi,m,k)-h_0
\mathbf{F}(\phi,k)\right),
$$
where
$$
\mathbf{F}(\phi,k)=\int\limits_0^\phi
\frac{dt}{\sqrt{1-k^2\sin^2t}}
$$
is the incomplete elliptic integral of the first kind,
$$
\mathbf{\Pi}(\phi,m,k)=\int\limits_0^\phi
\frac{dt}{(1+m\sin^2t)\sqrt{1-k^2\sin^2t}}
$$
is the incomplete elliptic integral of the third kind,
$$
l=C_1 \sqrt{\frac{c^*}{(p+q)(f(r)-F(\rho))}}\,,
$$
$$
\phi=\arcsin \frac{1}{k(1-c^*)}\,,\quad
m=k(c^*-1)\,,\quad h_0=\frac{f(r)}{h}.
$$

On the right-hand side
$$
L=\sqrt{B}\int\limits_\rho^\infty
\frac{d\zeta}{\sqrt{\zeta(\rho-\zeta)(1-\rho\zeta)}},
$$
changing the integration variable
$\tau=\rho/\zeta$, we have
$$
L=2\frac{\sqrt{B}}{\rho}
\mathbf{K}\left(\frac{1}{\rho}\right).
$$
Hence, integrating \eqref{eq3.4},
we obtain
\begin{equation}\label{eq3.6}
\left(l\mathbf{\Pi}(\phi,m,k)-h_0
\mathbf{F}(\phi,k)\right)=
\frac{\sqrt{B}}{\rho}
\mathbf{K}\left(\frac{1}{\rho}\right).
\end{equation}

Now in the $w$-plane, take the point $F(1)=f(-e^{i\alpha})$.
Integrate the equalities that are obtained from the system
\eqref{eq3.3}--\eqref{eq3.4} by extracting the square root the first over
$z$ from $0$ to $-1$, and then over the arc
$|z|=1$ counterclockwise from
$-1$ to $-e^{i\alpha}$, and the second,
over $\zeta$ from $1$ to $\rho$.

Write down the integrals on the left-hand side of \eqref{eq3.3}:
$$
\int\limits_0^{f(-1)}
\sqrt{\frac{C_1w-C_2}{w(w-F(\rho))(w-f(r))}}\,dw,
$$
$$
\int\limits_{f(-1)}^{F(1)}
\sqrt{\frac{C_1w-C_2}{w(w-F(\rho))(w-f(r))}}\,dw.
$$
Proceed with the integral on the left-hand side of \eqref{eq3.4}
$$
\int\limits_{F(1)}^{F(\rho)}
\sqrt{\frac{C_1w-C_2}{w(w-F(\rho))(w-f(r))}}\,dw.
$$
Summing up these integrals, we have
$$
T=\int\limits_0^{F(\rho)}
\sqrt{\frac{C_1w-C_2}{w(w-F(\rho))(w-f(r))}}\,dw.
$$
Making the change of variables $t=w/F(\rho)$, note that
$$
T=a^{*}\int\limits_0^1
\frac{t-b^{*}}{\sqrt{t(1-t)(1-c^{*}t)(t-b^{*})}}\,dt.
$$
Putting $t=1/x$, we find
$$
T=a^*\int\limits_1^\infty
\frac{dx}{x\sqrt{(x-1)(x-c^{*})(1-b^{*}x)}}-
a^*b^*\int\limits_1^\infty
\frac{dx}{\sqrt{(x-1)(x-c^{*})(1-b^{*}x)}}\,.
$$
Changing the integration variable
$y=b^{*}(x-1)/(b^{*}x-1)$ in the integrals and performing the corresponding transformations,
we obtain
$$
T=2\,i q h^{*}\mathbf{\Pi}(n^{*},k'),
$$
where
$$
h^*=\sqrt{\frac{(f(r)-F(\rho))F^2(\rho)}{C_2}}\,,
\quad k'=\sqrt{1-k^2},
\quad n^{*}=-\frac{1}{b^*}.
$$

Integrate the right-hand sides of \eqref{eq3.3}.
First, integrate over
$z$ from $0$ to $-1$
$$
T_{1}=-\sqrt{A}\int\limits_{-1}^0
\frac{dz}{\sqrt{z(r-z)(1-rz)}}.
$$
Performing the change of the integration variable
$x=(1+z)/(1-z)$ and applying a formula from
\cite{Gradshtein}, we obtain
$$
T_{1}=\sqrt{A}\,i\,
\mathbf{K}\left(\sqrt{1-r^{2}}\right).
$$

Integrate the right-hand side over the arc
$\beta$ of the unit circle counterclockwise from
$-1$ to $-e^{i\alpha}$:
$$
T_{2}=\sqrt{A}\int\limits_\beta
\frac{dz}{\sqrt{z(r-z)(1-rz)}}\;.
$$
Substituting $z=-e^{i\varphi}$ in $T_{2}$, where
$0\leq\varphi\leq \alpha$, we infer
$$
T_{2}=\sqrt{A}\int\limits_0^{\alpha}
\frac{d(-e^{i\varphi})}
{\sqrt{-e^{i\varphi}(r+e^{i\varphi})(1+re^{i\varphi})}}\;.
$$
Performing the change of variable
$t=-e^{i\varphi}$, we find
$$
T_{2}=-\frac{2\alpha}{\pi}\sqrt{A}\,
\mathbf{K}\left(r\right).
$$

Finally, integrate the right-hand side of \eqref{eq3.4} over $\zeta$ from $1$ to $\rho$:
$$
T_{3}=\sqrt{B}\int\limits_1^\rho
\frac{d\zeta}{\sqrt{\zeta(\rho-\zeta)(1-\rho\zeta)}}\;.
$$
Changing the integration variable
$u=\rho/\zeta$, we have
$$
T_{3}=\frac{\sqrt{B}}{\rho}\,i\,
\mathbf{K}\biggl(\sqrt{1-\frac{1}{\rho^{2}}}\biggr).
$$

Summing $T_{1}$, $T_{2}$ and $T_{3}$ yields
$$
\sqrt{A}\,i\,\mathbf{K}\left(\sqrt{1-r^{2}}\right)+
\frac{\sqrt{B}}{\rho}\,i\,
\mathbf{K}\biggl(\sqrt{1-\frac{1}{\rho^{2}}}\biggr)-
\frac{2\alpha}{\pi}\sqrt{A}\,\,
\mathbf{K}\left(r\right)\,.
$$

Thus, integrating the equalities
\eqref{eq3.3} and \eqref{eq3.4}, one has
$$
qh^{*}\mathbf{\Pi}(n^*,k')=\frac{\sqrt{A}}{2}\,
\mathbf{K}\left(\sqrt{1-r^{2}}\right)+
\\[3mm]
+
\frac{\sqrt{B}}{2\rho}\,
\mathbf{K}\biggl(\sqrt{1-\frac{1}{\rho^{2}}}\biggr)+
\frac{\alpha}{\pi}\,i
\sqrt{A}\,\mathbf{K}\left(r\right).
$$
Excluding the constants
$\sqrt{A}$ and $\sqrt{B}$ from this equality by using \eqref{eq3.5}
and \eqref{eq3.6}, we obtain the following equation
\begin{multline}\label{eq3.7}
-\frac{q h^*}{p h}\frac{ \mathbf{\Pi}(n^*,k')}{
\mathbf{\Pi}(n,k)}=\frac{1}{2}
\frac{\mathbf{K}(\sqrt{1-r^2})}{\mathbf{K}(r)}-
\\[3mm]
-\frac{1}{2}\,\frac{l\mathbf{\Pi}(\phi,m,k)-h_0
\mathbf{F}(\phi,k)}{p h \mathbf{\Pi}(n,k)}\,
\frac{\mathbf{K}\left(\sqrt{1-\frac{1}{\rho^2}}\right)}
{\mathbf{K}\left(\frac{1}{\rho}\right)}+
\frac{\alpha}{\pi}\,i.
\end{multline}

Hence, we have proved:

\begin{theorem}\label{Th.2}
The range $\Delta$ of the functional \eqref{eq2.1}
on the class $\mathfrak{M}$ is bounded by the curve defined by equation \eqref{eq3.7}.

\end{theorem}

\end{document}